# MIXED NONEUCLIDEAN GEOMETRIES

by Florentin Smarandache, Ph.D.
University of New Mexico
Department of Mathematics
Gallup, NM 87301, USA
Tel. (505) 726-1720
Fax: (505) 863-7532 (Attn. Smarandache)
E-mail: smarand@unm.edu

**Abstract**.
  The goal of this paper is to experiment new math concepts
and theories, especially if they run counter to the classical
ones.  To prove that contradiction is not a catastrophe, and
to learn to handle it in an (un)usual way.
To transform the apparently unscientific ideas into scientific
ones, and to develop their study (The Theory of Imperfections).
And finally, to interconnect opposite (and not only) human
fields of knowledge into as-heterogeneous-as-possible
another fields.

  The author welcomes any comments, notes, articles on this
paper and/or the 120 open questions bothering him,
which will be published in a collective monograph about
the paradoxist mathematics.

**Key words:**
    non-mathematics, anti-mathematics, dadaist algebra,
surrealist probability, cubist geometry, impressionist analysis,
theory of non-choice, wild algorithms, infinite computability
theory, symbolist mechanics, abstract physics, formalist
chemistry, expressionist statistics, hermetic combinatorics,
Sturm-und-Drang computer science, romanticist topology, letterist
number theory, illuminist set theory, aesthetic differential/
integral/functional equations, paradoxist logics, anti-literature,
experimental drama, non-poems, MULTI-STRUCTURE, MULTI-SPACE,
Euclidean spaces of non-integer or negative dimension, non-system,
anti-system, systems with infinetely many independent axioms,
unlimited theory, system of axioms based on a set with a single
element, INCONSISTENT SYSTEMS OF AXIOMS, CONTRADICTORY THEORY,
(unscientific, wrong, amalgam) geometry, (CHAOS or MESS)
GEOMETRIES (PARADOXIST GEOMETRY, NON-GEOMETRY, COUNTER-PROJECTIVE
GEOMETRY, ANTI-GEOMETRY), paradoxist model, critical area of a
model, paradoxist axioms, counter-axioms, counter-model,
counter-projective space, anti-axioms, anti-model, theory of
distorted buildings of Tits, paradoxist trigonometry,
DISCONTINUOUS MODELS, DISCONTINUOUS GEOMETRIES.

**1991 MSC:** 51D99

## INTRODUCTION:

These "Paradoxist Mathematics" may be understood as Experimental
Mathematics, Non-Mathematics, or even Anti-Mathematics:  not in a
nihilistic way, but in a positive one.  The truly innovative researchers
will banish the old concepts in order to check, by heuristic processes, some
new ones: their opposites.  Don't simply follow the crowd, and don't accept
to be manipulated by any (political, economical, social, even scientific,
or artistic, cultural, etc.) media!  Learn to contradict everything
and everybody!!
"Dubito, ergo cogito; cogito, ergo sum", said Rene Descartes, "I doubt,
 therefore I think; I think, therefore I exist" (metaphysical doubt).
See what happens if you deny the classics' theory?

   Since my childhood I didn't like the term of 'exact' sciences...
I hated it!  I didn't like the 'truth' displayed and given to me
on a plate -- as food to be swallowed although not to my taste.
   I considered the axioms as dogmas (not to think with your brain,
but with others'!), and I refused to follow them.
I wanted to be free in life -- because at that time I was
experiencing a political totalitarian system, without civil rights
-- hence I got the same feelings in science.
That's why I didn't trust anybody, especially the 'official' peoples.
(This is a REVOLT against all petrified knowledge.)
   A system of axioms means to me a dictatorship model in science.
It's not possible to perfectly formalize, i.e. without any intuition,
but sometimes researchers like to trick themselves!  Even Hilbert
recognized that just in his 1898 book of <Foundations of Geometry>
saying about the groups of axiomes that:  "Each of these groups
expresses, by itself, certain related fundamental facts of our
intuition".  And Kant in <Kritik der reinen Vernunft, Elementarlehre>,
Part 2, Sec. 2: "All human knowledge begins with intuition, thence
passes to concepts and ends with ideas".  Therefore, axiomatization
begins with intuition -- is it a paradox ?
The "traditional concept of recognizing the axioms as obvious truths
was replaced by the understanding that they are hypotheses for a
theory" [<Encyclopedic Dictionary of Mathematics>, second edition, by
the Mathematical Society of Japan, edited by Kiyosi Ito, translated
in English, MIT Press, Cambridge, Massachusetts, London, 1993, 35A,
p.155].
   The really avant-garde mind will entirely deny everything from
the past.  "No army can withstand the strenghth of an idea whose
time has come" (Victor Hugo).

## Questions 1-17 (one for each defined below section):
   While, in a usual way, people apply mathematics to other human
fields -- what about inserting literary and art theory in mathematics?
      How would we define the 'dadaist algebra', referring to the
         1916-22 nihilistic movement in literature, painting, sculpture
         that rejected all accepted conventions and produced non-sense
         and un-readable creations?  How can we introduce this style
         and similar <laws> in algebra??
      But the 'surrealist probability'? (this syntagme makes a little
         sense, doesn't it?).
      Or the 'cubist geometry', referring to the cubist paintings?
         (this may be exciting!).
      The 'impressionist analysis'?

```
        The 'theory of ... non-choice':
            - from two possibilities, pick the third one! (Buridan's ass!)
            - the best and unregrettable choice occurs when it's one and
                only possibility to choose from!
        The 'wild algorithms', meaning algorithms with an infinite number
            of (non-linear) steps;
        And the 'infinite computability theory' = how much of mathematics
            can be described in such wild algorithms.
        Same directions of study towards:
                'symbolist mechanics',
                'abstract physics' (suppose, for example, as an axiom,
                    that the speed of light is surpassed --
                    [see Homer B. Tilton, "Light beyond belief", Echo
                     Electronic Press, Tucson, 1995],
                    but if the speed of a material body can be unbounded,
                        even towards infinite?
                    and see what you get by this anti-relativity theory:
                    inventing new physics),
                'formalist chemistry',
                'expressionist statistics',
                'hermetic combinatorics',
                'Sturm-und-Drang computer science' (!)
                'romanticist topology' (wow, love is involving!)
                'letterist number theory' (!)
                'illuminist set theory',
                'esthetic differential/integral/functional equations',
        etc.
```

## Question 18:

```
        The 'paradoxist logics', referring to the F. Smarandache's 1980
            Paradoxist Literary Movement of avant-gardes, which may lead
            you to the anti-logic (which is logical!).
            Features of the 'paradoxist logics':
                # The Basic Thesis of the paradoxism:
                    everything has a meaning and a non-meaning
                    in a harmony each other.
                # The Essence of the paradoxism:
                        a) the sense has a non-sense,
                    and reciprocally
                        b) the non-sense has a sense.
                # The Motto of the paradoxism:
                    "All is possible, the impossible too!"
                # The Symbol of the paradoxism:
                    (a spiral -- optic illusion, or vicious circle)

                # The Delimitation from other avant-gardes:
                    - the paradoxism has a significance (in literature,
                    art, science),
                    while the dadaism, the lettrism, the absurd movement
                    do not;
                    - the paradoxism especially reveals the contradictions,
                    the anti-nomies, the anti-theses, the anti-phrases, the
                    antagonism, the non-conformism, the paradoxes in other
                    words of anything (in literature, art, science),
                    while the futurism, cubism, surrealism, abstractism
```

and all other avant-gardes do not focus on them.
           # The Directions of the paradoxism:
                - to use scientific methods (especially algorithms) for
                generating (and also studying) contradictory literary and
                artistic works;
                   and reciprocally
                to use artistic and literary methods for generating (and
                also studying) contradictory scientific works;
                - to create contradictory literary and artistic works in
                scientific spaces (using scientific: symbols, meta-language,
                matrices, theorems, lemmas, definitions, etc.).

## Question 19:

     From Anti-Mathematics to Anti-Literature:
     - I wrote a drama trilogy, called "MetaHistory", against the
       totalitarism of any kind: political, economical, social, cultural,
       artistic, even scientific (tendency of someones to monopolize the
       informational system, and to build not only political, economical,
       social dictatorships, but even dictatorships in culture, art, and
       science ... promoting only their people and friends, and
       boycotting the others);
          one of them, called "A Upside-Down World", with the property that
       by combinations of its scenes (which are independent modules) one
       gets 1,000,000,000 of bilions of different dramas!
          another drama, called "The Country of the Animals", has no ...
       dialogue!  (the characters' speech is showing on written placards).
     - I wrote "Non-poems":
             ^ poems with no words!
             ^ universal poems: poem-grafitti, poem-drawing, etc.
             ^ poems in 3-dimensional spaces
             ^ poems in Beltrami/Poincare/Hausdorff/etc. spaces
             ^ poetical models of ... mathematics:
                  poem-theorem, poem-lemma.
  Try the reverse way:  to apply math (and generally speaking science) in
arts and literarure.
(There are famous people:  as Lewis Carroll, Raymond Queneaux, Ion Barbu,
 etc. mathematicians and writers simultaneously.)

    Learn to deny (in a positive way) the masters and their work.  Thus
will progress our society.  Thus will make revolutionary steps towards
infinite ...
     Look at some famous examples:
- Lobacevsky contradicted Euclid in 1826:  "In geometry I find
certain imperfections", he said in his <Theory of Parallels>.
- Riemann came to contradict both his predecessors in 1854.
- Einstein contradicted Newton in early years of the XX-th century,
saying that if an object moves at a velocity close to the speed of
light, then time slows down, mass increases, and lenth in the direction
of motion decreases,
and so on...
Sometimes, people give new interpretations to old things ...
(and old interpretations to new things) !
[Don't talk about the humanistic field (art, literature, philosophy,
 sociology, etc.), where to reject other people's creation was and
 is being very common!  And much easier, comparing with the scientific
 field.]

What would be happened if everybody had obeyed the predecessors? (a stagnation).

**MULTI-STRUCTURE and MULTI-SPACES**.

I consider that life and practice do not deal with 'pure' spaces, but with a group of many spaces, with a mixture of structures, a 'mongrel', a heterogeneity -- the ardently preoccupation is to reunite them, to constitute a multi-structure.

I thought to a multi-space also: fragments (potsherds) of spaces put together, say as an example: Banach, Hausdorff, Tikhonov, compact, paracompact, Fock symmetric, Fock antisymmetric, path-connected, simply connected, discrete metric, indiscrete pseudo-metric, etc. spaces that work together as a whole mechanism. The difficulty is to be the passage over 'frontiers' (borders between two disjoint spaces); i.e. how can we organically tie a point P1 from a space S1 with a point P2 from a structurally opposite space S2 ? Does the problem become more complicated when the spaces' sets are not disjoint?

**Question 20**:
Can you define/construct Euclidian spaces of non-integer or negative dimension? [If so, are they connected in some way to Hausdorff's, or Kodaira's, Krull's, Lebesgue's (of a normal space), algebraic/cohomological(of a topological space, a scheme, or an associative algebra)/homological/(of a topological space, or a module) etc. dimension(s) ?]

**Question 21**:
Let's have the case of Euclid + Lobachevsky + Riemann geometric spaces (with corresponding structures) into one single space. What is the angles sum of a triangle with a vertex in each of these spaces equal to? and is it the same anytimes? Especialy to find a model of the below geometry would be interesting, or properties and applications of it. Paradoxically, the multi-, non-, or even anti- notions become after a while common notions. Their mystery, shock, novelty enter in the room of obvious things. This is the route of any invention and discovery.

Time is not uniform, but in a zigzag; a today's truth will be the tomorrow's falsehood -- and reciprocally, the opposite phenomena are complementary and may not survive independently. The every-day reality is a sumum or a multitude of rules, some of them opposite each other, accepted by ones and refused by others, on different surfaces of positive, negative, and null Gauss's curvatures in the same time (especially on non-constant curvature surfaces).

**Question 22**:
After all, what mathematical apparatus to use for subsequent improvement of this theory?

[my definition  is elementary].
        Logics without logics?
        System without system? (will be a non-system or anti-system?)
        Mathematics without mathematics!

   World is an ordered disorder and a disordered order!
Homogeneity exists only in pure sciences within our imagination,
but practice is quite different from theory.

   There are systems with one axiom only [see Dr. Paul Welsh, "Primitivity
in Mereology" (I and II), in <Notre Dame Journal of Formal Logic>, Vol.
XIX, No. 1 and 3, January and July 1978, pp. 25-62 and 355-85;
or B. Sobocinski, "A note on an axiom system of atomistic mereology", in
<Notre Dame Journal of Formal Logic>, Vol. XII, 1971, pp. 249-51.].
If one defines another system with a sole axiom, which is the negation of
the previous axiom, one gets an opposite theory.

## Question 23:
   Try to construct a consistent system of axioms, with infinitely
many independent axioms, in order to define a Unlimited Theory.
A theory to whom you may add at any time a new axiom to develop it in
all direction you like.

## Question 24:
    Try to construct a consistent system of axioms based on a set with
a single object (element).
(But if the set is ... empty?)

## INCONSISTENT SYSTEMS OF AXIOMS and CONTRADICTORY THEORY.

   Let (a1), (a2), ..., (an), (b) be n+1 independent axioms,
with n >= 1;  and let (b') be another axiom contradictory to (b).
We construct a system of n+2 axioms:
     [I ]      (a1), (a2), ..., (an), (b), (b')
which is inconsistent.  But this system may be shared into two
consistent systems of independent axioms
     [C ]      (a1), (a2), ..., (an), (b),
and
     [C']      (a1), (a2), ..., (an), (b').
We also consider the partial system of independent axioms
     [P ]      (a1), (a2), ... (an).
   Developing [P ], we find many propositions (theorems, lemmas)
              (p1), (p2), ..., (pm),
by combinations of its axioms.
   Developing [C ], we find all propositions of [P ]
              (p1), (p2), ..., (pm),
resulted by combinations of (a1), (a2), ..., (an),
plus other propositions
              (r1), (r2), ..., (rt),
resulted by combinations of (b) with any of (a1), (a2), ..., (an).
   Similarly for [C'], we find the propositions of [P ]
              (p1), (p2), ..., (pm),
plus other propositions
              (r'1), (r'2), ..., (r't),
resulted by combinations of (b') with any of (a1), (a2), ..., (an),

where (r'1) is an axiom contradictory to (r1), and so on.
  Now, developing [I ], we'll find all the previous resulted propositions:
            (p1), (p2), ..., (pm),
            (r1), (r2), ..., (rt),
            (r'1), (r'2), ..., (r't).
Therefore, [I ] is equivalent to [C ] reunited to [C'].
From one pair of contradictory propositions {(b) and (b')} in its
beginning, [I ] adds t more such pairs, where t >= 1, {(r1) and (r'1),
..., (rt) and (r't)}, after a complete step.
The further we go, the more pairs of contradictory propositions are
accumulating in [I ].

## Question 25:
  Develop the study of an inconsistent system of axioms.

## Question 26:
  It is interesting to study the case when n = 0.

Why do people avoid thinking about the **CONTRADICTORY THEORY** ?
As you know, nature is not perfect:
  and opposite phenomena occur together,
  and opposite ideas are simultaneously asserted and, ironically, proved
that both of them are true!  How is that possible? ...
A statement may be true in a referential system, but false in another
one.  The truth is subjective.  The proof is relative.
( In philosophy there is a theory:
  that "knowledge is relative to the mind, or things can be known only
  through their effects on the mind, and consequently there can be no
  knowledge of reality as it is in itself",
  called "the Relativity of Knowledge";
  <Webster's New World Dictionary of American English>, Third College
  Edition, Cleveland & New York, Simon & Schuster Inc., Editors:
  Victoria Neufeldt, David B. Guralnik, 1988, p. 1133. )
You know? ... sometimes is good to be wrong!

## Question 27:
  Try to develop a particular contradictory theory.

  I was attracted by Chaos Theory, deterministic behaviour which seems to
be randomly:  when initial conditions are varying little, the differential
equation solutions are varying tremendously much.  Originated by Poincare,
and studied by Lorenz, a metereologist, in 1963, by computer help.
These instabilities occuring in the numerical solutions of differential
equations are thus connected to the phenomena of chaos.
Look, I said, chaos in mathematics, like in life and world!
  Somehow consequently are the following four concepts in the paradoxist
mathematics, that may be altogether called, CHAOS (or MESS) GEOMETRIES!

**PARADOXIST GEOMETRY**

  In 1969, intrigued by geometry, I simultaneously constructed a
partially euclidean and partially non-euclidean space by a strange

replacement of the Euclid's fifth postulate (axiom of parallels)
with the following five-statement proposition:

   a) there are at least a straight line and a point exterior
      to it in this space for which only one line passes through
      the point and does not intersect the initial line;
      [1 parallel]
   b) there are at least a straight line and a point exterior
      to it in this space for which only a finite number of
      lines l , ..., l  (k >= 2) passe through the point and do not
            1        k
      intersect the initial line;
      [2 or more (in a finite number) parallels]
   c) there are at least a straight line and a point exterior
      to it in this space for which any line that passes through
      the point intersects the initial line;
      [0  parallels]
   d) there are at least a straight line and a point exterior
      to it in this space for which an infinite number of lines
      that passe through the point (but not all of them) do not
      intersect the initial line;
      [an infinite number of parallels, but not all lines passing
         throught]
   e) there are at least a straight line and a point exterior
      to it in this space for which any line that passes through
      the point does not intersect the initial line;
      [an infinite number of parallels, all lines throught
         the point]

I have called it the PARADOXIST GEOMETRY.
This geometry unites all together: Euclid, Lobachevsky/Bolyai, and
Riemann geometries.  And separates them as well!

Question 28:
  Now, the problem is to find a nice model (on manifolds) for this
Paradoxist Geometry, and study some of its characteristics.

## NON–GEOMETRY

  It's a lot easier to deny the Euclid's five postulates than
Hilbert's twenty thorough axioms.

   1. It is not always possible to draw a line from an arbitrary point
      to another arbitrary point.

      For example:
      this axiom can be denied only if the model's space has at least
      a discontinuity point;
      (in our bellow model MD, one takes an isolated point I in between
      f1 and f2, the only one which will not verify the axiom).

   2. It is not always possible to extend by continuity a finite line
      to an infinite line.

For example:
consider the bellow Model, and the segment AB, where both A and
B lie on f1, A in between P and N, while B on the left side of N;
one can not at all extend AB either beyond A or beyond B,
because the resulted curve, noted say A'-A-B-B', would not be a
geodesic (i.e. line in our Model) anymore.

If A and B lie in delta1-f1, both of them closer to f1, A in the
left side of P, while B in the right side of P,
then the segment AB, which is in fact A-P-B, can be extended
beyond A and also beyond B only up to f1
(therefore one gets a finite line too, A'-A-P-B-B', where A', B'
are the intersections of PA, PB respectively with f1).

If A, B lie in delta1-f1, far enough from f1 and P, such that AB
is parallel to f1, then AB verifies this postulate.

3. It is not always possible to draw a circle from an arbitrary
   point and of an arbitrary interval.

   For example:
   same as for the first axiom;
   the isolated point I, and a very small interval not reaching f1
   neither f2, will deny this axiom.

4. Not all the right angles are congruent.

   (See example of the Anti-Geometry, explained bellow.)

5. If a line, cutting two other lines, forms the interior angles of
   the same side of it strictly less than two right angles,
   then not always the two lines extended towards infinite cut each
   other in the side where the angles are strictly less than two right
   angles.

   For example:
   let h1, h2, l be three lines in delta1-delta2, where h1 intersects
   f1 in A, and h2 intersects f1 in B,  with A, B, P different each other,
   such that h1 and h2 do not intersect, but l cuts h1 and h2 and forms
   the interior angles of one of its side (towards f1) strictly less than
   two right angles;
   the assumption of the fifth postulate is fulfilled, but the
   consequence does not hold, because h1 and h2 do not cut each other
   (they may not be extended beyond A and B respectively, because the
   lines would not be geodesics anymore).

Question 29:
  Find a more convincing model for this non-geometry.

# COUNTER-PROJECTIVE GEOMETRY

Let P, L be two sets, and r a relation included in PxL.  The elements of
P are called points, and those of L lines.  When (p, l) belongs to r, we
say that the line l contains the point p.
For these, one imposes the following COUNTER-AXIOMS:

(I) There exist:   either at least two lines, or no line,
      that contains two given distinct points.

(II) Let p1, p2, p3 be three non-collinear points,  and q1, q2 two
      distinct points.  Suppose that {p1, q1, p3} and {p2, q2, p3} are
      collinear triples.  Then the line containing p1, p2, and the line
      containing q1, q2 do not intersect.

(III) Every line contains at most two distinct points.

## Questions 30-31:
Find a model for the Counter-(General Projective) Geometry
(the previous I and II counter-axioms hold),
  and a model for the Counter-Projective Geometry
(the previous I, II, and III counter-axioms hold).
[They are called **COUNTER-MODELS** for the general projective, and projective
geometry, respectively.]

## Questions 32-33:
Find geometric models for each of the following two cases:
  - There are points/lines that verify all the previous counter-axioms,
      and other points/lines in the same **COUNTER-PROJECTIVE SPACE**
      that do not verify any of them;
  - Some of the counter-axioms I, II, III are verified,
      while the others are not
      (there are particular cases already known).

## Question 34:
The study of these counter-models may be extended to Infinite-Dimensional
Real (or Complex) Projective Spaces, denying the IV-th axiom, i.e.:
  (IV) There exists no set of finite number of points for which any
      subspace that contains all of them contains P.

## Question 35:
Does the Duality Principle hold in a counter-projective space?
  What about Desargues's Theorem, Fundamental Theorem of Projective
Geometry / Theorem of Pappus, and Staudt Algebra ?
  Or Pascal's Theorem, Brianchon's Theorem ?
(I think none of them will hold!)

## Question 36:
The Theory of Buildings of Tits, which contains the Projective
Geometry as a particular case, can be 'distorted' in the same <paradoxist>
way by deforming its axiom of a BN-pair (or Tits system) for the triple
(G, B, N), where G is a group, and B, N its subgroups;
[see J. Tits, "Buildings of spherical type and finite BN-pairs", Lecture notes in math. 386, Springer, 1974].
Notions as: simplex, complex, chamber, codimension, apartment,
building will get contorted either ...
   Develop a Theory of Distorted Buildings of Tits!

## ANTI-GEOMETRY

   It is possible to entirely de-formalize Hilbert's groups of axioms
of the Euclidean Geometry, and to construct a model such that none of
his fixed axioms holds.
      Let's consider the following things:
         - a set of <points>:  A, B, C, ...
         - a set of <lines>:  h, k, l, ...
         - a set of <planes>:  alpha, beta, gamma, ...
      and
         - a set of relationships among these elements:  "are situated",
            "between", "parallel", "congruent", "continuous", etc.
Then, we can deny all Hilbert's twenty axioms [see David Hilbert,
"Foundations of Geometry", translated by E. J. Townsend, 1950;
and Roberto Bonola, "Non-Euclidean Geometry", 1938].
There exist cases, within a geometric model, when the same axiom is
verified by certain points/lines/planes and denied by others.

      GROUP I.  ANTI-AXIOMS OF CONNECTION:

         I.1.  Two distinct points A and B do not always completely
               determine a line.

               Let's consider the following model MD:
               get an ordinary plane delta, but with an infinite
               hole in of the following shape:

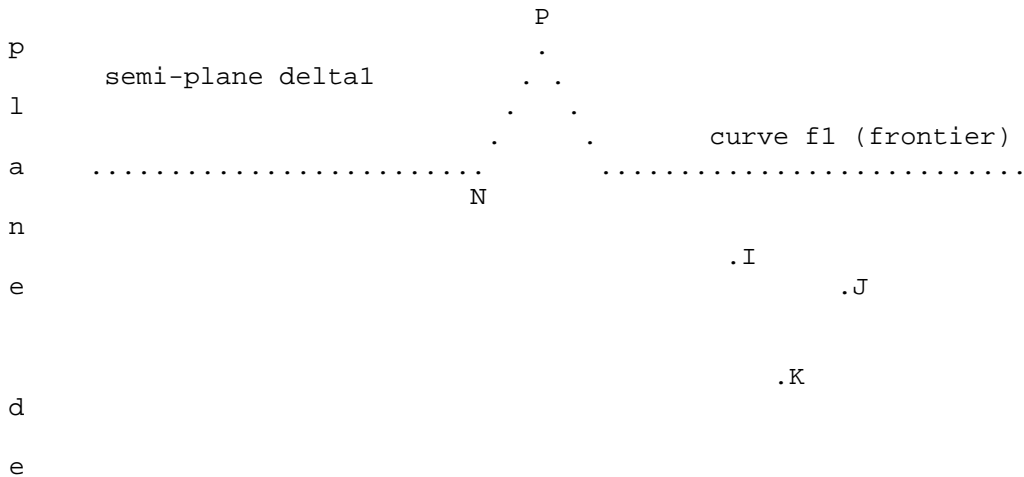

```
                                        P
  p                                     .
        semi-plane delta1             . .
  l                                  .   .
                                    .      .        curve f1 (frontier)
  a   ........................           ............................
                                   N
  n
                                         .I
  e                                              .J

                                            .K
  d

  e
```

```
l
     .........................        ............................
t                              .    .        curve f2 (frontier)
                                .  .
a      semi-plane delta2          . .
                                  .
                                  Q
```

            Plane delta is a reunion of two disjoint planar
            semi-planes;
            f1 lies in MD, but f2 does not;
            P, Q are two extreme points on f that belong to MD.

               One defines a LINE l as a geodesic curve:  if two
            points A, B that belong to MD lie in l, then
            the shortest curve lied in MD between A and B
            lies in l also.
            If a line passes two times through the same point, then
            it is called double point (KNOT).

               One defines a PLANE alpha as a surface such that for
            any two points A, B that lie in alpha and belong to
            MD there is a geodesic which passes through A, B and
            lies in alpha also.

            Now, let's have two strings of the same length:
            one ties P and Q with the first string s1 such that
            the curve s1 is folded in two or more different
            planes and s1 is under the plane delta;
            next, do the same with string s2, tie Q with P, but
            over the plane delta and such that s2 has a different
            form from s1;
            and a third string s3, from P to Q, much longer than s1.
            s1, s2, s3 belong to MD.

            Let I, J, K be three isolated points -- as some islands,
            i.e. not joined with any other point of MD,
            exterior to the plane delta.

               This model has a measure, because the (pseudo-)line is
            the shortest way (lenth) to go from a point to another
            (when possible).

**Question 37**:
            Of course, this model is not perfect, and is far from
            the best.  Readers are asked to improve it, or to make up
            a new one that is better.

            (Let A, B be two distinct points in delta1-f1.  P and Q are
            two points on s1, but they do not completely determine a
            line, referring to the first axiom of Hilbert,
            because A-P-s1-Q are different from B-P-s1-Q.)

        I.2.  There is at least a line l and at least two distinct

```
            points A and B of l, such that A and B do not
            completely determine the line l.

            (Line A-P-s1-Q are not completely determined by P and Q
             in the previous construction, because B-P-s1-Q is another
             line passing through P and Q too.)

I.3.    Three points A, B, C not situated in the same line do
        not always completely determine a plane alpha.

        (Let A, B be two distinct points in delta1-f1, such that
         A, B, P are not co-linear.  There are many planes
         containing these three points: delta1 extended with any
         surface s containing s1, but not cutting s2 in between
         P and Q, for example.)

I.4.    There is at least a plane, alpha, and at least three
        points A, B, C in it not lying in the same line, such
        that A, B, C do not completely determine the plane
        alpha.

        (See the previous example.)

I.5.    If two points A, B of a line l lie in a plane alpha,
        doesn't mean that every point of l lies in alpha.

        (Let A be a point in delta1-f1, and B another point on
         s1 in between P and Q.  Let alpha be the following plane:
         delta1 extended with a surface s containing s1, but not
         cutting s2 in between P and Q, and tangent to delta2 on
         a line QC, where C is a point in delta2-f2.
         Let D be point in delta2-f2, not lying on the line QC.
         Now, A, B, D are lying on the same line A-P-s1-Q-D,
         A, B are in the plane alpha, but D do not.)

I.6.    If two planes alpha, beta have a point A in common,
        doesn't mean they have at least a second point in
        common.

        (Construct the following plane alpha: a closed surface
         containing s1 and s2, and intersecting delta1 in one point
         only, P.  Then alpha and delta1 have a single point in
         common.)

I.7.    There exist lines where lies only one point,
        or planes where lie only two points,
        or space where lie only three points.

        (Hilbert's I.7 axiom may be contradicted if the model
         has discontinuities.
         Let's consider the isolated points area.
           The point I may be regarded as a line, because it's not
```

```
           possible to add any new point to I to form a line.
              One constructs a surface that intersects the model only
           in the points I and J.)

GROUP II.  ANTI-AXIOMS OF ORDER:

  II.1.  If A, B, C are points of a line and B lies between A
         and C, doesn't mean that always B lies also between
         C and A.

         [Let T lie in s1, and V lie in s2, both of them
          closer to Q, but different from it.  Then:
            P, T, V are points on the line P-s1-Q-s2-P
            ( i.e. the closed curve that starts from the point P
              and lies in s1 and passes through the point Q and
              lies back to s2 and ends in P ),
            and T lies between P and V
              -- because PT and TV are both geodesics --,
            but T doesn't lie between V and P
              -- because from V the line goes to P and then to T,
                 therefore P lies between V and T.]

         [By definition: a segment AB is a system of points
          lying upon a line between A and B (the extremes are
          included).

          Warning:  AB may be different from BA;
          for example:
            the segment PQ formed by the system of points
            starting with P, ending with Q, and lying in s1,
            is different from the segment QP formed by the
            system of points starting with Q, ending with P,
            but belonging to s2.
          Worse, AB may be sometimes different from AB;
          for example:
            the segment PQ formed by the system of points
            starting with P, ending with Q, and lying in s1,
            is different from the segment PQ formed by the
            system of points starting with P, ending with Q,
            but belonging to s2.]

  II.2.  If A and C are two points of a line, then:
           there does not always exist a point B lying between A
            and C,
           or there does not always exist a point D such that C lies
            between A and D.

         [For example:
          let F be a point on f1, F different from P,
          and G a point in delta1, G doesn't belong to f1;
          draw the line l which passes through G and F;
          then:
           there exists a point B lying between G and F
             -- because GF is an obvious segment --,
```

```
                but there is no point D such that F lies between
                G and D -- because GF is right bounded in F
                ( GF may not be extended to the other side of F,
                 because otherwise the line will not remain a
                 geodesic anymore ).]

II.3.   There exist at least three points situated on
        a line such that:
             one point lies between the other two,
             and another point lies also between the other two.

        [For example:
         let R, T be two distinct points, different
         from P and Q, situated on the line P-s1-Q-s2-P,
         such that the lenghts PR, RT, TP are all equal;
         then:
            R lies between P and T,
            and T lies between R and P;
            also P lies between T and R.]

II.4.   Four points A, B, C, D of a line can not always be
        arranged:
              such that B lies between A and C and also
               between A and D,
              and such that C lies between A and D and also between
               B and D.

        [For examples:
         - let R, T be two distinct points, different
         from P and Q, situated on the line P-s1-Q-s2-P such
         that the lenghts PR, RQ, QT, TP are all equal,
         therefore R belongs to s1, and T belongs to s2;
         then P, R, Q, T are situated on the same line:
            such that R lies between P and Q, but not between
             P and T
               -- because the geodesic PT does not pass through
                  R --,
            and such that Q does not lie between P and T
               -- because the geodesic PT does not pass through
                  Q --,
             but lies between R and T;
         - let A, B be two points in delta2-f2 such that A, Q, B
         are colinear, and C, D two points on s1, s2 respectively,
         all of the four points being different from P and Q;
         then A, B, C, D are points situated on the same line
         A-Q-s1-P-s2-Q-B, which is the same with line
         A-Q-s2-P-s1-Q-B, therefore we may have two different
         orders of these four points in the same time:
         A, C, D, B and A, D, C, B.]

II.5.   Let A, B, C be three points not lying in the same
        line, and l a line lying in the same plane ABC and
        not passing through any of the points A, B, C.
        Then, if the line l passes through a point of the
```

segment AB, it doesn't mean that always the line l
will pass through either a point of the segment BC
or a point of the segment AC.

[For example:
 let AB be a segment passing through P in the
 semi-plane delta1, and C a point lying in delta1
 too on the left side of the line AB;
 thus A, B, C do not lie on the same line;
 now, consider the line Q-s2-P-s1-Q-D, where D is
 a point lying in the semi-plane delta2 not on f2:
 therefore this line passes through the point P of
 the segment AB, but do not pass through any point
 of the segment BC, nor through any point of the
 segment AC.]

GROUP III.   ANTI-AXIOM OF PARALLELS.

        In a plane alpha there can be drawn through a point
        A, lying outside of a line l, either no line,
        or only one line, or a finite number of lines,
        or an infinite number of lines which do not intersect
        the line l.  (At least two of these situations should occur.)
        The line(s) is (are) called the parallel(s) to l
        through the given point A.

        [ For examples:
         - let l0 be the line N-P-s1-Q-R, where N is a point
            lying in delta1 not on f1, and R is a similar
            point lying in delta2 not on f2,
            and let A be a point lying on s2,
            then:  no parallel to l0 can be drawn through A
            (because any line passing through A, hence through
            s2, will intersect s1, hence l0, in P and Q);
         - if the line l1 lies in delta1 such that l1 does
            not intersect the frontier f1, then:
            through any point lying on the left side of l1
            one and only one parallel will pass;
         - let B be a point lying in f1, different from P,
            and another point C lying in delta1, not on f1;
            let A be a point lying in delta1 outside of BC;
            then:  an infinite number of parallels to the
            line BC can be drawn through the point A.

          Theorem.  There are at least two lines l1, l2 of a
            plane, which do not meet a third line l3 of the
            same plane, but they meet each other,
            ( i.e. if l1 is parallel to l3, and l2 is parallel
              to l3, and all of them are in the same plane,
              it's not necessary that l1 is parallel to l2 ).
            [ For example:
              consider three points A, B, C lying in f1, and
              different from P, and D a point in delta1 not on
              f1;  draw the lines AD, BE and CE such that
              E is a point in delta1 not on f1 and both BE

```
                and CE do not intersect AD;
                then:  BE is parallel to AD, CE is also parallel
                to AD, but BE is not parallel to CE because the
                point E belong to both of them. ]

GROUP IV.   ANTI-AXIOMS OF CONGRUENCE

        IV.1. If A, B are two points on a line l, and A' is a
              point upon the same or another line l', then:
              upon a given side of A' on the line l', we can
              not always find only one point B' so that the
              segment AB is congruent to the segment A'B'.

              [ For examples:
              - let AB be segment lying in delta1 and having no
                point in common with f1, and construct the line
                C-P-s1-Q-s2-P (noted by l') which is the same
                with C-P-s2-Q-s1-P, where C is a point lying in
                delta1 not on f1 nor on AB;
                take a point A' on l', in between C and P, such
                that A'P is smaller than AB;
                now, there exist two distinct points B1' on s1
                and B2'on s2, such that A'B1' is congruent to AB
                and A'B2' is congruent to AB,
                with A'B1' different from A'B2';
              - but if we consider a line l' lying in delta1 and
                limited by the frontier f1 on the right side
                (the limit point being noted by M),
                and take a point A' on l', close to M, such that
                A'M is less than A'B', then:  there is no point
                B' on the right side of l' so that A'B' is
                congruent to AB. ]

              A segment may not be congruent to itself!

              [ For example:
                - let A be a point on s1, closer to P,
                  and B a point on s2, closer to P also;
                  A and B are lying on the same line A-Q-B-P-A
                  which is the same with line A-P-B-Q-A,
                  but AB meseared on the first representation
                  of the line is strictly greater than AB
                  meseared on the second representation of
                  their line. ]

        IV.2.  If a segment AB is congruent to the segment
               A'B' and also to the segment A''B'', then
               not always the segment A'B' is congruent to
               the segment A''B''.

               [ For example:
               - let AB be a segment lying in delta1-f1, and
                 consider the line C-P-s1-Q-s2-P-D, where C, D are
                 two distinct points in delta1-f1 such that C, P, D
                 are colinear.  Suppose that the segment AB is
                 congruent to the segment CD (i.e. C-P-s1-Q-s2-P-D).
```

Get also an obvious segment A'B' in delta1-f1,
different from the preceding ones, but congruent
to AB.
Then the segment A'B' is not congruent to the segment
CD (considered as C-P-D, i.e. not passing through Q.)

IV.3.   If AB, BC are two segments of the same line l
        which have no points in common aside from the
        point B,
        and A'B', B'C' are two segments of the same line
        or of another line l' having no point other than
        B' in common, such that AB is congruent to A'B'
        and BC is congruent to B'C',
        then not always the segment AC is congruent to
        A'C'.

     [ For example:
        let l be a line lying in delta1, not on f1,
        and A, B, C three distinct points on l, such
        that AC is greater than s1;
        let l' be the following line: A'-P-s1-Q-s2-P
        where A' lies in delta1, not on f1,
        and get B' on s1 such that A'B' is congruent
        to AB, get C' on s2 such that BC is congruent
        to B'C' (the points A, B, C are thus chosen);
        then:  the segment A'C' which is first seen as
        A'-P-B'-Q-C' is not congruent to AC,
        because A'C' is the geodesic A'-P-C' (the
        shortest way from A' to C' does not pass
        through B') which is strictly less than AC. ]

Definitions.    Let h, k be two lines having a point O
                in common.  Then the system (h, O, k) is
                called the angle of the lines h and k in
                the point O.
                ( Because some of our lines are curves,
                  we take the angle of the tangents to
                  the curves in their common point. )

                The angle formed by the lines h and k
                situated in the same plane, noted by
                <(h, k), is equal to the arithmetic mean
                of the angles formed by h and k in all
                their common points.

IV.4.   Let an angle (h, k) be given in the plane alpha,
        and let a line h' be given in the plane beta.
        Suppose that in the plane beta a definite side
        of the line h' be assigned, and a point O'.
        Then in the plane beta there are one, or more,
        or even no half-line(s) k' emanating from the
        point O' such that the angle (h, k) is
        congruent to the angle (h', k'),
        and at the same time the interior points of

the angle (h', k') lie upon one or both sides
of h'.

[ Examples:
 - Let A be a point in delta1-f1, and B, C two
   distinct points in delta2-f2;
   let h be the line A-P-s1-Q-B, and k be the
   line A-P-s2-Q-C;  because h and k intersect
   in an infinite number of points (the segment
   AP), where they normally coincide -- i.e. in
   each such point their angle is congruent to
   zero, the angle (h, k) is congruent to zero.
   Now, let A' be a point in delta1-f1, different
   from A, and B' a point in delta2-f2, different
   from B, and draw the line h' as A'-P-s1-Q-B';
   there exist an infinite number of lines k', of
   the form A'-P-s2-Q-C' (where C' is any point in
   delta2-f2, not on the line QB'), such that the
   angle (h, k) is congruent to (h', K'),
   because (h', k') is also congruent to zero, and
   the line A'-P-s2-Q-C' is different from the line
   A'-P-s2-Q-D' if D' is not on the line QC'.
 - If h, k, and h' are three lines in delta1-P,
   which intersect the frontier f1 in at most one
   point, then there exist only one line k' on a
   given part of h' such that the angle (h, k) is
   congruent to the angle (h', k').
 - *Is there any case when, with these hypotheses,
   no k' exists ?
 - Not every angle is congruent to itself;
   for example:
      <(s1, s2) is not congruent to <(s1, s2)
   [because one can construct two distinct lines:
    P-s1-Q-A and P-s2-Q-A, where A is a point in
    delta2-f2, for the first angle, which becomes equal
    to zero;
    and P-s1-Q-A and P-s2-Q-B, where B is another point
    in delta2-f2, B different from A, for the second
    angle, which becomes strictly greater than zero!].

IV. 5. If the angle (h, k) is congruent to the angle
       (h', k',) and the angle (h'', k''), then the
       angle (h', k') is not always congruent to the
       angle (h'', k'').

       (A similar construction to the previous one.)

IV. 6. Let ABC and A'B'C' be two triangles such that
          AB is congruent to A'B',
          AC is congruent to A'C',
          <BAC is congruent to <B'A'C'.
       Then not always
          <ABC is congruent to <A'B'C'
          and <ACB is congruent to <A'C'B'.

[For example:
       Let M, N be two distinct points in delta2-f2, thus
       obtaining the triangle PMN;
       Now take three points R, M', N' in delta1-f1, such
       that RM' is congruent to PM, RN' is congruent to RN,
       and the angle (RM', RN') is congruent to the angle
       (PM, PN).  RM'N' is an obvious triangle.
       Of course, the two triangle are not congruent,
       because for example PM and PN cut each other twice
       -- in P and Q -- while RM' and RN' only once -- in
       R.
       (These are geodesical triangles.)]

Definitions:

  Two angles are called supplementary if they have the
  same vertex, one side in common, and the other sides
  not common form a line.

  A right angle is an angle congruent to its
  supplementary angle.

  Two triangles are congruent if its angles are congruent
  two by two, and its sides are congruent two by two.

Propositions:

  A right angle is not always congruent to another
  right angle.

  For example:
  Let A-P-s1-Q be a line, with A lying in delta1-f1,
  and B-P-s1-Q another line, with B lying in
  delta1-f1 and B not lying in the line AP;
  we consider the tangent t at s1 in P, and B chosen
  in a way that <(AP, t) is not congruent to <(BP, t);
  let A', B' be other points lying in delta1-f1
  such that <APA' is congruent to <A'P-s1-Q,
  and <BPB' is congruent to <B'P-s1-Q.
  Then:
   - the angle APA' is right, because it is congruent
   to its supplementary (by construction);
   - the angle BPB' is also right, because it is
   congruent to its supplementary (by construction);
   - but <APA' is not congruent to <BPB',
   because the first one is half of the angle A-P-s1-Q,
   i.e. half of <(AP, t),
   while the second one is half of the B-P-s1-Q,
   i.e. half of <(BP, t).

  The theorems of congruence for triangles
  [side, side, and angle in between;  angle, angle, and
  common side;  side, side, side] may not hold either
  in the Critical Zone (s1, s2, f1, f2) of the Model.

Property:
The sum of the angles of a triangle can be:
- 180 degrees, if all its vertexes A, B, C are
lying, for example, in delta1-f1;
- strictly less than 180 degrees [ any value in the
interval (0, 180) ],
for example:
let R, T be two points in delta2-f2 such that Q does
not lie in RT, and S another point on s2;
then the triangle SRT has <(SR, ST) congruent to 0
because SR and ST have an infinite number of common
points (the segment SQ), and <QTR + <TRQ congruent
to 180-<TQR [ by construction we may vary <TQR in the
interval (0, 180) ];
- even 0 degree!
let A be a point in delta1-f1, B a point in delta2-f2,
and C a point on s3, very close to P;
then ABC is a non-degenerate triangle (because its
vertexes are non-colinear), but <(A-P-s1-Q-B, A-P-s3-C)
= <(B-Q-s1-P-s3-C) = <(C-s3-P-A,
C-s3-P-s1-Q-B) = 0
(one considers the lenth C-s3-P-s1-Q-B strictly less
than C-s3-B);
the area of this triangle is also 0 !
- more than 180 degrees,
for example:
let A, B be two points in delta1-f1, such that
<PAB + <PBA + <(s1, s2; in Q) is strictly greater
than 180 degrees;
then the triangle ABQ, formed by the intersection of
the lines A-P-s2-Q, Q-s1-P-B, AB will have the sum of
its angles strictly greater than 180 degrees.

Definition:
A circle of center M is a totality of all points A
for which the segments MA are congruent to one another.

For example, if the center is Q, and the lenth of the
segments MA is chosen greater than the lenth of s1,
then the circle is formed by the arc of circle centered
in Q, of radius MA, and lying in delta2, plus another
arc of circle centered in P, of radius MA-lenth of s1,
lying in delta1.

GROUP V. ANTI-AXIOM OF CONTINUITY (ANTI-ARCHIMEDEAN AXIOM)

Let A, B be two points. Take the points A1, A2, A3, A4,
... so that A1 lies between A and A2, A2 lies between
A1 and A3, A3 lies between A2 and A4, etc. and the
segments AA1, A1A2, A2A3, A3A4, ... are congruent to one
another.
Then, among this series of points, not always  there exists
a certain point An such that B lies between A and An.

For example:
let A be a point in delta1-f1, and B a point on f1, B

different from P;
on the line AB consider the points A1, A2, A3, A4, ...
in between A and B, such that AA1, A1A2, A2A3, A3A4, etc.
are congruent to one another;
then we find that there is no point behind B (considering
the direction from A to B), because B is a limit point
(the line AB ends in B).

The Bolzano's (intermediate value) theorem may not hold in
the Critical Zone of the Model.

## Question 38:

It's very interesting to find out if this system of axioms are complete
and consistent (!)  The apparent unscientific or wrong geometry, which looks
more like an amalgam, is somehow supported by its attached model.

## Question 39:

  How will the differential equations look like in this field?

## Question 40:

  How will the (so called by us:) "PARADOXIST" TRIGONOMETRY look like in this
field?

## Question 41:

  First, one can generalize this model using more bridges (conections/
strings between delta1 and delta2) of many lenths, and many gates (points
like P and Q on f1 and f2, respectively) -- from a finite to an infinite
number of such bridges and gates.
If one put all bridges in the delta plane, one gets a dimension-2 model;
otherwise, the dimension is >= 3.
Some bridges may be replaced with (round or not necessaryly) bodies,
tangent (or not necesaryly) to the frontiers f1 and f2.

## Question 42:

  Should it be indicated to remove the discontinuities ?
  But what about DISCONTINUOUS MODELS (on spaces not everywhere
continuous -- like our MD) ? generating in this way
DISCONTINUOUS GEOMETRIES.

## Question 43:

  The model MD can also be generalized to n-dimensional space as a
hypersurface, considering the group of all projective transformations of
an (n+1)-dimensional real projective space that leave MD invariant.

## Question 44-47:

  Find geometric models for each of the following four cases:
    - NO point/line/plane in the model space verifies any of Hilbert's
       twenty axioms;
       (in our MD, some points/lines/planes did verify, and some others
       did not);
    - The Hilbert's groups of axioms I, II, IV, V are denied for any point/
       line/plane in the model space, but the III-th one (axiom of
       parallels) is verified;
       this is an Opposite-(Lobacevski+Riemann) Geometry:
       neither hyperbolic, nor elliptic ... and yet Non-Euclidean!

- The groups of anti-axioms I, II, IV, V are all verified, but the
    III-th one (anti-axiom of parallels) is denied;
- Some of the groups of anti-axioms I, II, III, IV, V are verified,
    while the others are not -- except the previous case;
    (there are particular cases already known).

## Question 48:

  What connections may be found among this Paradoxist Model,
and the Cayley, Klein, Poincare, Beltrami (differential geometric)
models?

## Question 49-120 (combinig by twos, each new geometry -- out of 4 -- with an old geometry -- out of 18 -- all mentionned below):

  What connections among these Paradoxist Geometry, Non-Geometry, Counter-
Projective Geometry, Anti-Geometry  and the other ones:  Conformal (Mobius)
Geometry, Pseudo-Conformal Geometry, Laguerre Geometry, Desarguesian and
Non-Desarguesian Geometries, Non-Archimedian Geometry, Spectral Geometry,
Spherical Geometry, Hyper-Sphere Geometry, Wave Geometry (Y. Mimura),
Non-Holonomic Geometry (G. Vranceanu), Cartan's Geometry of Connection,
Integral Geometry (W. Blaschke), Continuous Geometry (von Neumann),
Affine Geometry, Generalized Geometries (of H. Weyl, O. Veblen, J. A.
Schoutten), etc.

## CONCLUSION:

  The above 120 OPEN QUESTIONS are not impossible at all.
"The world is moving so fast nowadays that the person, who says <it
can't be done>, is often interrupted by someone doing it"!
[ <Leadership> journal, Editor Arthur F. Lenehan, October 24, 1995, p. 16,
  Fairfield, NJ ].

  The author encourages readers to send not only comments, but also new
(solved or unsolved) questions arising from them.

  Specials thanks to professors JoAnne Growney, Zahira S. Khan, and Paul
Hartung of Bloomsburg University, Pennsylvania, for giving me the
opportunity to write this article and to lecture it on November 13th, 1995,
in their Department of Mathematics and Computer Sciences.